# BAYESIAN VARIABLE SELECTION FOR HIGH DIMENSIONAL GENERALIZED LINEAR MODELS: CONVERGENCE RATES OF THE FITTED DENSITIES

By Wenxin Jiang

*Northwestern University*

Bayesian variable selection has gained much empirical success recently in a variety of applications when the number $K$ of explanatory variables $(x_1, \ldots, x_K)$ is possibly much larger than the sample size $n$. For generalized linear models, if most of the $x_j$'s have very small effects on the response $y$, we show that it is possible to use Bayesian variable selection to reduce overfitting caused by the curse of dimensionality $K \gg n$. In this approach a suitable prior can be used to choose a few out of the many $x_j$'s to model $y$, so that the posterior will propose probability densities $p$ that are "often close" to the true density $p^*$ in some sense. The closeness can be described by a Hellinger distance between $p$ and $p^*$ that scales at a power very close to $n^{-1/2}$, which is the "finite-dimensional rate" corresponding to a low-dimensional situation. These findings extend some recent work of Jiang [Technical Report 05-02 (2005) Dept. Statistics, Northwestern Univ.] on consistency of Bayesian variable selection for binary classification.

**1. Introduction.** Bayesian variable selection (BVS) is a fruitful method for studying regression models that relate a response $y$ to a vector of candidate explanatory variables $x = (x_1, \ldots, x_K)^T$. For example, when generalized linear models (GLM) are considered, the density of $y$ and the mean function of $y$ conditional on $x$ both depend on a linear combination $x^T \beta$ through the regression coefficients $\beta = (\beta_1, \ldots, \beta_K)^T$. The BVS approach uses priors that propose different model $\gamma$'s and the corresponding sets of regression coefficient $\beta_\gamma$'s, where $\gamma$ indicates the components of $x$ that are included in regression. The posterior distribution $\pi[\gamma, \beta_\gamma | D^n]$ for the model and the model parameters $(\gamma, \beta_\gamma)$ can then be obtained based on an observed data set $D^n = (x^{(i)}, y^{(i)})_1^n$, which is often assumed to consist of i.i.d. (independent and identically distributed) copies of $(x, y)$. Computational simplification is achievable in the cases of linear regression and probit regression, where the unknown regression coefficients $\beta_\gamma$ can often be analytically integrated out in the posterior-based computations (e.g., Kohn, Smith and Chan [17]; Lee et al. [19]).

The BVS approach has had many successful applications. For example, when applied in a linear regression framework, BVS is used in basis selection for nonparametric regression (e.g., Smith and Kohn [23], Kohn, Smith and Chan [17]) and in construction of financial index tracking portfolios (e.g., George and McCulloch [7]). Other work applying BVS in the GLM framework includes, for example, Clyde and DeSimone-Sasinowska [3], Nott and Leonte [21] and Wang and George [24]. Recently, BVS has been applied to the area of bioinformatics. In order to construct Gaussian graphical models for gene expression pathways, Dobra et al. [5] obtain biologically meaningful results by applying Bayesian variable selection to model how each gene in the graph relates to tens of thousands of other genes. In order to classify binary responses based on microarray data, Lee et al. [19] and Sha et al. [22] (via probit regression) and Zhou, Liu and Wong [27] (via logistic regression) use BVS to achieve excellent cross-validated classification errors. These most recent applications are especially noteworthy since they are all in the situation of $K \gg n$, where the number of candidate variables $K$ can be several thousand and the sample size $n$ is often less than a hundred.

Despite these empirical successes, there has not been a systematic study of the frequentist properties of BVS, such as posterior consistency and convergence rates. It is the aim of this paper to study these convergence properties for BVS, allowing $K$ to be possibly much larger than $n$. The consistency that we will consider is neither the traditional sense (i) of consistency in estimating the true regression parameters, nor the sense (ii) of consistency in identifying the true model (the $x$-components with nonzero regression coefficients). Sense (i) is not feasible since in cases with $K \gg n$, the $\beta$-coefficients are often not identifiable. Sense (ii) is not a totally satisfactory framework when, as in many realistic situations, none of the $K$ regression coefficients is exactly zero, even though many of them may be very small. The consistency we consider is the closeness between the true (conditional) density $p^* = p^*(y|x)$ and the densities $p = p(y|x; \gamma, \beta_\gamma)$ proposed by the posterior $\pi(\gamma, \beta_\gamma | D^n)$. We do not attempt to identify the "true parameter" or the "true model" (the nonzero coefficients). Rather, we allow all coefficients to be not exactly zero, and attempt to construct the posterior to propose models that include only a few of those nonzero coefficients, but have the corresponding densities $p$ "often close" to the true $p^*$ in some sense.

Let $\nu_x(dx)$ be the probability measure for $x$ and $\nu_y(dy)$ be the dominating measure for conditional densities $p$ and $p^*$. Define the Hellinger distance between $p$ and $p^*$ as $d(p,p^*) = \sqrt{\int\int \nu_y(dy)\nu_x(dx)(\sqrt{p}-\sqrt{p^*})^2}$. The convergence results we consider describe the "often closeness" between $p$ and $p^*$ that can be formulated as, for example,

(1) $$P^*[\pi[d(p,p^*) < \varepsilon_n | D^n] > 1 - \delta_n] \geq 1 - \lambda_n,$$

for all large enough $n$, for some small $\varepsilon_n, \delta_n, \lambda_n$ converging to zero as $n \to \infty$, where $P^*$ is the probability measure for data $D^n$, when they are generated as i.i.d. copies from $p^*\nu_y(dy)\nu_x(dx)$.

It is noted that BVS is essential for achieving convergence results as above, when $K \gg n$. A usual approach without variable selection, using the full model and putting a prior on all the regression coefficients, can be shown to lead to bad results in the following counterexample.

EXAMPLE. Suppose $K > n$. Let the random variable $z$ take values from $\{j/K\}_{j=1}^n$ with equal probability and let $x = (x_1,\ldots,x_K)^T$, where $x_j = I[z = j/K]$ for each $j$. Let $(z^{(i)}, x^{(i)}, y^{(i)})_{i=1}^n$ and $(z,x,y)$ be i.i.d., where $y|x \sim N(0,1)$. Suppose the fitted model is $y|x \sim N(\sum_{j=1}^K \beta_j x_j, 1)$, where, without selecting among the $(x_j)_1^K$, one proposes a prior for $(\beta_j)_1^K$ as i.i.d. $N(0,1)$.

The Hellinger distance $d(p,p^*)$ between $p^* = e^{-y^2/2}/\sqrt{2\pi}$ and $p = e^{-(y-\sum_{j=1}^K \beta_j x_j)^2/2}/\sqrt{2\pi}$ is such that $d(p,p^*)^2 = (2/K)\sum_{j=1}^K (1 - e^{-\beta_j^2/8})$.

Then, in the posterior conditional on $D^n = (x^{(i)}, y^{(i)})_{i=1}^n$, $\beta_1,\ldots,\beta_K$ are independent and $\beta_j \sim N(\sum_{i=1}^n y^{(i)} x_j^{(i)}/(1+\sum_{i=1}^n x_j^{(i)}), 1/(1+\sum_{i=1}^n x_j^{(i)}))$. Note that $(x_j^{(i)})_{i=1}^n$ are zero for at least $K-n$ of the $K$ $j$'s since $x_j^{(i)} = I[z^{(i)} = j/K]$, and the $n$ $z^{(i)}$'s can only populate at most $n$ out of $K$ of the $j/K$ locations. Therefore, at least $K-n$ out of the $K$ $\beta_j$'s follow the $N(0,1)$ distribution in the posterior—which is the same as the corresponding prior distribution. Without loss of generality, let $\beta_1^{K-n}$ be independent $N(0,1)$ in the posterior. Note that $d(p,p^*)^2 \geq [2(K-n)/K][1/(K-n)]\sum_1^{K-n}(1-e^{-\beta_j^2/8})$. For a simple example, let $K = 2n$. An application of Chebyshev's inequality leads to $\pi[d(p,p^*) \geq \eta^{1/2}|D^n] \geq 1 - (\eta^2 n)^{-1}$, for $\eta = 1/2 - 1/\sqrt{5}$, which happens with $P^*$-probability 1. Therefore, without variable selection, a convergence result such as (1) cannot hold. Such a convergence result, however, can be shown to hold for this example, with $\varepsilon_n$ following a near finite-dimensional rate (a power close to $1/\sqrt{n}$), if Bayesian variable selection is used properly, according to later results of this paper (e.g., Remark 1).

There has been considerable interest recently in studying the theoretical properties of high-dimensional regression. Most results are for frequentist methods. For example, Bühlmann [2] considers boosting for high-dimensional regression; Greenshtein and Ritov [12] and Greenshtein [11] consider constrained or $\ell_1$-penalized optimization; Meinshausen and Bühlmann [20] apply a similar method of $\ell_1$ penalization to high-dimensional graphical models. Recently, Fan and Li [6] have provided a useful overview for methods based on the penalized likelihood for treating high dimensionality, which includes examples of generalized linear models and survival models, among others. In contrast to these frequentist approaches based on optimization, the Bayesian method considered here has the attractive capability of presenting several likely models together with the corresponding posterior probabilities. A theoretical study of Bayesian inference *without* variable selection has been carried out by, for example, Ghosal [8, 9]. This work considers $K$'s growing with $n$ but at a slower rate. On the other hand, in the $K \ll n$ case treated in Ghosal [8, 9], posterior asymptotic normality was established for the whole parameter vector, so the goal there was much higher, and hence, the result there is not comparable with the result in the present paper which focuses on posterior convergence rates.

In contrast to previous work, we consider Bayesian variable selection and allow the cases $K \gg n$. It is noted that it is essential to have the variable selection step in order to obtain good results when $K \gg n$. The counterexample above shows that, without variable selection, it is impossible to have good convergence in general cases with $K > n$, while with variable selection, excellent empirical performance has been reported, for example, in Lee et al. [19] and Sha et al. [22] with $K \gg n$.

We study the convergence behavior of BVS for generalized linear models, which include linear regression, logistic regression, probit regression, Poisson regression, and so on. We also include a discussion of Gaussian graphical models that uses linear regression for neighborhood selection. Therefore, the current paper forms an extension to Jiang [15], who only considers consistency of BVS for binary logistic and probit regression, without studying the convergence rates. Here we study the convergence rate $\varepsilon_n$ as well, and will show that despite the high-dimension $K \gg n$, BVS can still lead to a near finite-dimensional rate (with $\varepsilon_n$ close to $1/\sqrt{n}$ in order), if we are in some "sparse" situations when most of the regression coefficients are very small. (For binary regression, this rate $\varepsilon_n$ also forms a good convergence rate for the purpose of classification, as shown in Section 5 later.) For such sparse high-dimensional problems, Bayesian variable selection can therefore help to reduce "overfitting" or the "curse of dimensionality." Note that such a conclusion can only be drawn by a careful study of the convergence rates in high dimensions; just proving the consistency, as in Jiang [15], is not enough; for example, it is well known (e.g., Hastie, Tibshirani and Friedman [13], Chapter 13) that the $k$-nearest neighbor rules are consistent for classification, but can suffer considerably from the curse of dimensionality. Also, it is well known (e.g., Devroye, Györfi and Lugosi [4], Chapters 6 and 7) that (at least in finite dimensions) there exist universally consistent classification rules, but any rule can have a very slow convergence rate under some data distribution.

Below we will first specify the notation and the framework of the paper.

**2. Notation and framework.** The explanatory variable is a $K_n$-dimensional random vector $x = (x_1, x_2, \ldots, x_{K_n})^T$. Following the typical practice of studying high-dimensional problems, we will formally consider the asymptotics when $K_n$ increases as $n \to \infty$.

For simplicity, we will assume that $|x_j| \leq 1$ for all $j$ for most of the later discussion. The results can be easily extended to the case when all $|x_j|$'s are bounded above by a large constant.

The response is $y$. The true relation between $y$ and $x$ is assumed to follow a parametric generalized linear model (GLM) with true conditional density $p^*(y|x)$ and the corresponding mean function $\mu^*(x)$. Generalized linear models (GLM) are a class of popular regression models relating a response $y$ to a vector of covariates $x$. The GLM with one natural parameter is constructed with a density of the form $p^*(y|x) = \exp\{a(h^*)y + b(h^*) + c(y)\} \equiv f(y, h^*)$, where $h^* = x^T \beta^*$ is the linear parameter, $a(h)$ and $b(h)$ are continuously differentiable, and $a(h)$ has nonzero derivative. The mean function $\mu^* = E(y|x) = -b'(h^*)/a'(h^*) \equiv \psi(x^T \beta^*)$ follows a transformed linear model, where the transform $\psi$ is the inverse of a chosen link function. This formalism includes regression models for responses that are binary, Poisson and Gaussian (with known error variance), and can be easily extended to the cases with a dispersion parameter, which can then include Gaussian models with unknown error variance.

We assume that corresponding to the true model $p^*$, there exists a true regression parameter vector $\beta^*$, which satisfies some "sparseness" conditions, describing situations when most components of $\beta^*$ are very small in magnitude. One such condition states that $\lim_{n \to \infty} \sum_{j=1}^{K_n} |\beta_j^*| < \infty$. Other conditions can be formulated to describe how fast the sum of $|\beta_j^*|$'s converges.

The condition limiting the sum of $|\beta_j^*|$'s has been considered by Bühlmann [2] for studying how boosting algorithms handle high-dimensional linear regression. As Bühlmann points out, as a special case, this condition is satisfied when only a finite and fixed number of $x_j$'s are relevant, that is, when the number of nonzero $\beta_j^*$'s is independent of $n$. More generally, the sparseness conditions can describe situations when all $x_j$'s are relevant, but most of them have very small effects ($|\beta_j^*|$'s).

Note that $p^*$ is the conditional density of $y|x$, which is also the joint density of $(x, y)$ if the dominating measure $\nu_x(dx)\nu_y(dy)$ is the product of the probability measure of $x$ and the dominating measure of $y$. We will always use this kind of dominating measure.

The data for $n$ subjects are assumed to be independent and identically distributed (i.i.d.) based on $p^*\nu_x(dx)\nu_y(dy)$. Therefore, showing the subject index $i$, the data set is of the form $D^n = (x_1^{(i)}, \ldots, x_{K_n}^{(i)}, y^{(i)})_{i=1}^n$.

The prior selects a subset of the $K_n$ $x$-variables in the data set to model $y$, using density $f(y, x_\gamma^T \beta_\gamma)$, where $\gamma = (\gamma_1, \ldots, \gamma_{K_n})$ has 0/1 valued components which are 1 only when the corresponding $x$ component is included in the model, that is, $\gamma_j = I[|\beta_j| > 0]$. Sometimes we will also use $\gamma$ to denote the corresponding set of index $j$'s for which $|\beta_j| > 0$. The notation $v_\gamma$ denotes the subvector of a vector $v$ with components $\{v_j\}$, for all $j$'s with $\gamma_j = 1$ (or for all $j \in \gamma$, if $\gamma$ is understood as the corresponding index set).

We use the probability measure $\pi_n(\gamma, d\beta_\gamma)$ to denote the prior distribution of the subset model $\gamma$ and the corresponding regression coefficients $\beta_\gamma$. (The prior depends on the sample size $n$, but we will often drop the subscript $n$ for simpler notation.) This induces a posterior measure conditional on the data set $D^n$,

$$\pi(\gamma, d\beta_\gamma | D^n)$$
$$= \prod_{i=1}^n p(y_i, x_i | \gamma, \beta_\gamma) \pi(\gamma, d\beta_\gamma) \Big/ \sum_{\gamma'} \int_{\beta'_{\gamma'}} \prod_{i=1}^n p(y_i, x_i | \gamma', \beta'_{\gamma'}) \pi(\gamma', d\beta'_{\gamma'}),$$

where $p(y, x | \gamma, \beta_\gamma) = f(y, x_\gamma^T \beta_\gamma)$. The prior and posterior distributions for $(\gamma, \beta_\gamma)$ induce distributions for the corresponding parameterized densities.

For notational simplification, we will use $|v|$ to denote the sum of the absolute values of the components for any vector $v$. For two positive sequences $a_n$ and $b_n$, $a_n \prec b_n$ (or $b_n \succ a_n$) means $\lim_{n \to \infty} a_n / b_n = 0$.

**3. A prior specification.** General conditions on the prior will be given later in Section 7. Here, for being specific, we first consider the following prior for $(\gamma, \beta_\gamma)$. Conditional on $\gamma$, $\beta_\gamma$ follows $N(0, V_\gamma)$, where $V_\gamma$ is a $|\gamma| \times |\gamma|$ covariance matrix.

To complete this prior specification, we let the model indicators $\gamma = (\gamma_1, \ldots, \gamma_{K_n})$ be generated by first proposing i.i.d. binary random variables $\tilde{\gamma}_1^n$, with $\pi(\tilde{\gamma}_j = 1) = \lambda_n = r_n / K_n$, where we assume, for convenience, that $r_n$ is some integer smaller than $K_n$. We then keep only the $\tilde{\gamma}$'s satisfying a size restriction $\sum |\tilde{\gamma}_j| \leq \bar{r}_n$, and let the prior proposed model $\gamma = \tilde{\gamma}$. Here $r_n$ is the prior expectation of model size $|\tilde{\gamma}|$ *before* applying the size restriction; $\bar{r}_n$ is the maximal possible model size. We assume $1 \leq r_n \leq \bar{r}_n < K_n$.

Therefore, $\pi(\gamma) \propto \prod_{j=1}^{K_n} \lambda_n^{\gamma_j} (1 - \lambda_n)^{1 - \gamma_j} I[\sum_{l=1}^{K_n} \gamma_l \leq \bar{r}_n]$. Although the size restriction is not necessary (see more general conditions in Section 7), it helps to keep the model from becoming too complicated and gives a convenient starting point for proving the theoretical properties. Also, without this kind of restriction, the design matrix $\sum_{i=1}^n x_{i\gamma} x_{i\gamma}^T$ would become singular when the proposed model size $|\gamma| > n$; such a design matrix is often used in the popular algorithms for generating the posterior distributions in Gaussian regression (e.g., Smith and Kohn [23]), probit regression (e.g., Lee et al. [19]) and logistic regression (e.g., Zhou et al. [27]).

Under this specification of prior, we will present conditions on $V_\gamma$, $r_n$ and $\bar{r}_n$ for proving results on posterior consistency and convergence rates. The condition on $V_\gamma$ will depend on how the largest eigenvalues ($\mathrm{ch}_1$) of $V_\gamma$ and $V_\gamma^{-1}$ grow with the size of $|\gamma|$.

Let $H(\gamma) = \max\{\mathrm{ch}_1(V_\gamma), \mathrm{ch}_1(V_\gamma^{-1})\}$. In many typical cases $H(\gamma)$ grows at most polynomially in model size, that is, $H(\gamma) \leq B|\gamma|^v$ for some constant $B > 0$ and some power $v > 0$, for all large $|\gamma|$. For example, when $V_\gamma = cI_\gamma$ (proportional to the identity matrix; see, e.g., Dobra et al. [5]), $H(\gamma)$ is a constant and does not grow with $|\gamma|$. For another example, $V_\gamma \approx constant \times (Ex_\gamma x_\gamma^T)^{-1}$; see, for example, Smith and Kohn [23] and Lee et al. [19], who use a sample approximation of this choice. Then the largest eigenvalues of $V_\gamma$ and $V_\gamma^{-1}$ are both bounded linearly for large $|\gamma|$, when $x_\gamma$ has components standardized to have mean zero and common variance, and have all pairwise correlations being $\rho \in (0, 1)$. In addition, for $V_\gamma$ following the covariance matrix of a finite-order AR or MA process, when the lag polynomials have no zeros on the unit circle, the eigenvalues of $V_\gamma$ and $V_\gamma^{-1}$ are also bounded such that $\max\{\mathrm{ch}_1(V_\gamma), \mathrm{ch}_1(V_\gamma^{-1})\}$ grows like $|\gamma|^0$. For a detailed discussion of these eigenvalues, see, for example, Section 3 of Bickel and Levina [1].

**4. Convergence results for GLM.** Here, for simplicity we will assume that all explanatory variables are bounded and standardized such that $|x_j| \leq 1$ for all $j$. Assume $\lim_{n \to \infty} \sum_1^{K_n} |\beta_j^*| < \infty$ for a regression parameter $\beta^*$ corresponding to the true density $p^*$, where $K_n$ is a nondecreasing sequence in $n$.

We also assume that the prior specification in Section 3 is used. Define $\Delta(r_n) = \inf_{\gamma:|\gamma|=r_n} \sum_{j:j \notin \gamma} |\beta_j^*|$, $B(r_n) = \sup_{\gamma:|\gamma|=r_n} \mathrm{ch}_1(V_\gamma^{-1})$ and $\bar{B}(r_n) = \sup_{\gamma:|\gamma|=r_n} \mathrm{ch}_1(V_\gamma)$. Let $\tilde{B}_n = \sup_{\gamma:|\gamma| \leq \bar{r}_n} \mathrm{ch}_1(V_\gamma)$. Let $D(R) = 1 + R \times \sup_{|h| \leq R} |a'(h)| \cdot \sup_{|h| \leq R} |\psi(h)|$ for any $R > 0$.

THEOREM 1. *Assume that the prior specification in Section 3 is used, $|x_j| \leq 1$ for all $j$ and $\lim_{n \to \infty} \sum_1^{K_n} |\beta_j^*| < \infty$, where $K_n$ is a nondecreasing sequence in $n$.*

*Let $\varepsilon_n$ be a sequence such that $\varepsilon_n \in (0, 1]$ for each $n$ and $n\varepsilon_n^2 \succ 1$ and assume that the following conditions also hold:*

(2) $$\bar{r}_n \ln(1/\varepsilon_n^2) \prec n\varepsilon_n^2,$$

(3) $$\bar{r}_n \ln(K_n) \prec n\varepsilon_n^2,$$

(4) $$\bar{r}_n \ln D(\bar{r}_n \sqrt{n\varepsilon_n^2 \tilde{B}_n}) \prec n\varepsilon_n^2,$$

(5) $$1 \leq r_n \leq \bar{r}_n < K_n,$$

(6) $$1 \prec r_n \prec K_n,$$

(7) $$\Delta(r_n) \prec \varepsilon_n^2,$$

(8) $$B(r_n) \prec n\varepsilon_n^2,$$

(9) $$r_n \ln \bar{B}(r_n) \prec n\varepsilon_n^2.$$

Denote $d(p,p^*)^2 = \int \int |p(y,x|\gamma,\beta_\gamma)^{1/2} - p^*(y,x)^{1/2}|^2 \nu_y(dy)\nu_x(dx)$. Then we have the following successively stronger results:

(i) for some $r_0 > 0$,
$$\lim_{n\to\infty} P^*\{\pi[d(p,p^*) \leq \varepsilon_n|D^n] \geq 1 - e^{-r_0 n\varepsilon_n^2}\} = 1;$$

(ii) for some $c_1 > 0$, and for all sufficiently large $n$,
$$P^*\{\pi[d(p,p^*) > \varepsilon_n|D^n] \geq e^{-0.5c_1 n\varepsilon_n^2}\} \leq e^{-0.5c_1 n\varepsilon_n^2};$$

(iii) for some $c_1 > 0$, and for all sufficiently large $n$,
$$E^*_{D^n}\pi[d(p,p^*) > \varepsilon_n|D^n] \leq e^{-c_1 n\varepsilon_n^2}.$$

The above condition on $D(\cdot) = D(\bar{r}_n\sqrt{n\varepsilon_n^2 \tilde{B}_n})$, when considering a specific example of GLM, depends on how $|a'(h)|$ and $|\psi'(h)|$ grow with the linear parameter $h$. We will consider the following examples here.

(a) *Poisson regression* with log linear link: mean $\mu = e^h$, $y \in \{0,1,2,\ldots\}$. Then
$$f(y,h) = \frac{e^{-\mu}}{y!}\mu^y = \exp\{hy - e^h - \ln(y!)\}.$$

Here $a(h) = h$, $a' = 1$, $\psi(h) = e^h$. So both $|a'|$ and $|\psi|$ grow at most exponentially in $|h|$.

(b) *Normal linear regression*: mean $\mu = h$; variance $\sigma^2 = \varphi^{-1} \in \Re^+$ is assumed to be known for now; $y \in \Re$. Then
$$f(y,h) = \frac{1}{\sqrt{2\pi\sigma^2}}e^{(-1/(2\sigma^2))(y-\mu)^2}$$
$$= \exp\left\{\varphi hy - \frac{\varphi h^2}{2} - \frac{\varphi y^2}{2} - \frac{1}{2}\ln(2\pi\varphi^{-1})\right\}.$$

Here $a(h) = \varphi h$, $a' = \varphi$, $\psi(h) = h$. So both $|a'|$ and $|\psi|$ grow at most linearly in $|h|$.

(c) *Exponential regression* with log linear link: mean $\mu = e^h$, $y \in (0,\infty)$. Then
$$f(y,h) = \mu^{-1}e^{-y/\mu} = \exp\{-e^{-h}y - h\}.$$

Here $a(h) = -e^{-h}$, $a' = e^{-h}$, $\psi = e^h$. So both $|a'|$ and $|\psi|$ grow at most exponentially in $|h|$.

(d) *Binary logistic regression*: mean $\mu = e^h/(1+e^h)$, $y \in \{0,1\}$. Then

$$f(y,h) = \mu^y(1-\mu)^{1-y} = \exp\{hy - \ln(1+e^h)\}.$$

Here $a(h) = h$, $a' = 1$, $\psi(h) = e^h/(1+e^h)$. So both $|a'|$ and $|\psi|$ are bounded above by 1.

(e) *Binary probit regression*: mean $\mu = \Phi(h) \equiv \int_{-\infty}^{h}(e^{-z^2/2}/\sqrt{2\pi})\,dz$, $y \in \{0,1\}$. Then

$$f(y,h) = \mu^y(1-\mu)^{1-y} = \exp\{y\ln(\Phi(h)/(1-\Phi(h))) + \ln(1-\Phi(h))\}.$$

Here $a(h) = \ln(\Phi(h)/(1-\Phi(h)))$, $a' = [\Phi(h)^{-1} + \{1-\Phi(h)\}^{-1}]\Phi'(h)$, $\psi(h) = \Phi(h) \in [0,1]$. By using Mills' ratio, it can be shown that $|a'(h)|$ increases at most linearly with $|h|$.

Using these rates of growth, we can make the condition on $D(\cdot)$ more specific for specific examples of GLM.

The conditions of Theorem 1 also depend on how the eigenvalues of $V_\gamma$ and $V_\gamma^{-1}$ behave. To be specific, assume that the largest eigenvalues of $V_\gamma$ and $V_\gamma^{-1}$, for $|\gamma| \leq \bar{r}_n$, are both bounded above by some power $\bar{r}_n^v$ ($v > 0$), for all large enough $\bar{r}_n$.

The condition on $r_n \ln \bar{B}(r_n)$ then becomes redundant since $r_n \ln \bar{B}(r_n) \leq \bar{r}_n \ln \tilde{B}_n \leq c\bar{r}_n \ln \bar{r}_n \leq c\bar{r}_n \ln K_n \prec n\varepsilon_n^2$ (for some constant $c > 0$ and for all large enough $n$), since $\tilde{B}_n$ is bounded above by a power of $\bar{r}_n$.

The condition on $\bar{r}_n \ln(1/\varepsilon_n^2)$ also becomes redundant [they are implied by the condition on $\bar{r}_n \ln(K_n)$ and $n\varepsilon_n^2 \succ 1$] if we assume that $K_n \succ n^\delta$ for some $\delta > 0$.

Consider now the condition on $\bar{r}_n \ln D(\cdot)$ for various regression models, depending on the rate of growth $D(\bar{r}_n\sqrt{n\varepsilon_n^2 \tilde{B}_n})$. This condition on $\bar{r}_n \ln D(\cdot)$ becomes redundant [it is implied by the condition on $\bar{r}_n \ln(K_n)$] for normal linear regression, binary logistic regression, and probit regression since $D(\cdot)$ is bounded above by some power of $\sqrt{\bar{r}_n^{2+v} n\varepsilon_n^2}$, which is bounded above by some power of $K_n$ (note that $\bar{r}_n \leq K_n$, $\varepsilon_n \leq 1$ and $K_n \succ n^\delta$ for some $\delta > 0$). The condition $B(r_n) \prec n\varepsilon_n^2$ can be satisfied by requiring $\bar{r}_n \prec (n\varepsilon_n^2)^{1/v}$.

For Poisson and exponential regressions with the log-linear link, however, $D(\cdot)$ grows exponentially in $\sqrt{\bar{r}_n^{2+v} n\varepsilon_n^2}$. The condition on $\bar{r}_n \ln D(\cdot)$ then cannot be ignored, and it can be satisfied if $\bar{r}_n \prec (n\varepsilon_n^2)^{1/(4+v)}$ [which actually implies the later condition $B(r_n) \prec n\varepsilon_n^2$ and makes it redundant].

These are summarized as follows.

THEOREM 2. *Assume that the prior specification in Section 3 is used, such that $\max\{\sup_{\gamma:|\gamma|\leq\bar{r}_n} \mathrm{ch}_1(V_\gamma), \sup_{\gamma:|\gamma|\leq\bar{r}_n} \mathrm{ch}_1(V_\gamma^{-1})\} \leq B\bar{r}_n^v$ for some positive constants $B$ and $v$, for all large enough $\bar{r}_n$. Suppose $|x_j| \leq 1$ for all $j$,*

and $\lim_{n\to\infty} \sum_1^{K_n} |\beta_j^*| < \infty$, where $K_n$ is a nondecreasing sequence in $n$ and $K_n \succ n^\delta$ for some $\delta > 0$.

Let $\varepsilon_n$ be a sequence such that $\varepsilon_n \in (0,1]$ for each $n$ and $n\varepsilon_n^2 \succ 1$ and assume that the following conditions also hold:

$$\bar{r}_n \ln(K_n) \prec n\varepsilon_n^2, \tag{10}$$

$$1 \leq r_n \leq \bar{r}_n < K_n, \tag{11}$$

$$1 \prec r_n \prec K_n, \tag{12}$$

$$\Delta(r_n) \equiv \inf_{\gamma:|\gamma|=r_n} \sum_{j:j\notin\gamma} |\beta_j^*| \prec \varepsilon_n^2. \tag{13}$$

Also assume that

$$\bar{r}_n \prec (n\varepsilon_n^2)^{1/v} \tag{14}$$

for normal linear, binary logistic and binary probit regression; or assume

$$\bar{r}_n \prec (n\varepsilon_n^2)^{1/(4+v)} \tag{15}$$

for Poisson or exponential regression with log-linear link function. Then the results of Theorem 1 hold.

This result can be used to study the convergence rate $\varepsilon_n$ under various situations, depending on how $K_n$ grows with $n$, as well as how $\Delta(r_n) = \inf_{|\gamma|=r_n} \sum_{j\notin\gamma} |\beta_j^*|$ grows with $r_n$. Here are some corollaries, which follow by assuming an exponential decay rate of $\Delta(\cdot)$ and checking the conditions of Theorem 2. This includes as a special case only a fixed and finite number of $|\beta_j^*|$'s being nonzero, while also allowing a more realistic setup with many small $|\beta_j^*|$'s, none of which is exactly zero.

COROLLARY 1. *Consider the examples of Poisson regression, exponential regression, normal linear regression, logistic regression or probit regression described before. Assume that the prior specification in Section 3 is used, such that $\max\{\sup_{\gamma:|\gamma|\leq \bar{r}_n} \text{ch}_1(V_\gamma), \sup_{\gamma:|\gamma|\leq \bar{r}_n} \text{ch}_1(V_\gamma^{-1})\} \leq B\bar{r}_n^v$ for some positive constants $B$ and $v$, for all large enough $\bar{r}_n$. Suppose $|x_j| \leq 1$ for all $j$. Suppose $K_n \succ n^\delta$ for some $\delta > 0$ and $K_n \leq e^{Cn^\xi}$ for some $C > 0$ and some $\xi \in (0,1)$, for all large enough $n$. Suppose $\lim_{n\to\infty} \sum_{j=1}^{K_n} |\beta_j^*| < \infty$. Also suppose for some $C' > 0$, $\Delta(r_n) \leq e^{-C'r_n}$ for all large enough $n$, and*

$$(C')^{-1} \ln n \leq r_n \leq \bar{r}_n \prec (\ln n)^k \tag{16}$$

*for some $k > 1$. Then we can take the convergence rate in Theorem 2 as*

$$\varepsilon_n \sim n^{-(1-\xi)/2} (\ln n)^{k/2}. \tag{17}$$

REMARK 1 (*Good convergence rate*). Note that $n^\alpha \prec e^{Cn^\xi}$ for any small $\xi > 0$ and large $\alpha > 0$. So if $K_n \sim n^\alpha$ for whatever large power $\alpha$, one can achieve a convergence rate $\varepsilon_n \sim n^{-(1-\xi)/2}(\ln n)^{k/2} \prec n^{-(1-2\xi)/2}$, where $\xi$ can be made arbitrarily close to zero. This gives a rate arbitrarily close to the "finite-dimensional" rate $1/\sqrt{n}$, despite the large dimension $K_n$. We note also that these results suggest slowly growing $r$ and $\bar{r}_n$ between powers of $\ln n$, for achieving a near finite-dimensional convergence rate. Since these are only sufficient conditions, it may be possible that other ranges of $r$ and $\bar{r}_n$ can also lead to a near finite-dimensional rate of convergence. The following result, for example, shows a good convergence rate even when $r$ and $\bar{r}_n$ grow slowly in some small power of $n$.

COROLLARY 2. *Consider the setup of Corollary* 1. *For any $b \in (0, q)$, if instead of (16) we have*

$$(18) \qquad (C')^{-1} \ln n \leq r_n \leq \bar{r}_n \prec n^b,$$

*then we can take the convergence rate as*

$$(19) \qquad \varepsilon_n \sim n^{-(1-\xi-b)/2}.$$

*Here the power $q = \min\{1 - \xi, \delta, \xi/(3 + v)\}$ for Poisson and exponential regressions with log-linear link function; $q = \min\{1 - \xi, \delta\}I[v \leq 1] + \min\{1 - \xi, \delta, \xi/(v-1)\}I[v > 1]$ for logistic, probit and normal linear regression.*

REMARK 2 (*Posterior consistency*). The results on posterior consistency can be obtained as a special case by setting $\varepsilon_n = \varepsilon$ for any small but fixed $\varepsilon > 0$. There is no need to assume a rate for $\Delta(r_n)$ for consistency results to hold, since $\Delta(r_n) \prec \varepsilon^2$ as long as $r_n \succ 1$ and $\lim_{n\to\infty} \sum_1^{K_n} |\beta_j^*| < \infty$. The previous Theorem 2 then implies that the following condition on $r_n$ and $\bar{r}_n$ is sufficient for posterior consistency:

$$(20) \qquad 1 \prec r_n \leq \bar{r}_n \prec \min\{K_n, n^{1/(v+4)}, n/(\ln K_n)\}.$$

A slightly more relaxed condition for consistency for the special cases of logistic and probit regression can be found in Jiang [15].

REMARK 3 (*Normal linear regression with unknown dispersion*). So far, for normal linear regression, we have assumed that $y|x \sim N(E(y|x), \varphi^{-1})$ with dispersion parameter (inverse variance) $\varphi(>0)$ known. In practice, $\varphi$ is unknown and a gamma prior is often put on $\varphi$ (e.g., George and McCulloch [7], Kohn, Smith and Chan [17] and Dobra et al. [5]). For example, suppose conditional on model $\gamma$, $\varphi|\gamma \sim Ga(\kappa, \rho)$ with prior density $\pi(\varphi|\gamma) = \rho^\kappa \varphi^{\kappa-1} e^{-\rho\varphi}/\Gamma(\kappa)$, $\beta_\gamma|\gamma, \varphi \sim N(0, \varphi^{-1}V_\gamma)$ and $\gamma$ follows the prior distribution of Section 3. With this prior specification, it can be shown that the statements regarding normal linear regression in Theorem 2 and Corollaries 1 and 2 are still valid, where we consider bounded covariates standardized such that $|x_j| \leq 1$ for all $j$.

**5. Implications of posterior convergence.** It is well-known that a convergence statement such as

$$\lim_{n \to \infty} P^*\{\pi_n[d(p, p^*) \leq \varepsilon_n | D^n] \geq 1 - e^{-r_0 n \varepsilon_n^2}\} = 1, \tag{21}$$

for some $r_0 > 0$, implies existence of point estimates of $p^*$ that have the same convergence rate $\varepsilon_n$ in the frequentist sense. Such a point estimate can be obtained by finding the center of an $\varepsilon_n$-ball with high posterior probability, or by posterior expectation (e.g., Ghosal, Ghosh and van der Vaart [10]).

A point estimate can also be formed by a generalization of posterior expectation called a "selected posterior estimate" (Jiang [15]). For example,

$$\hat{p}_A = \int p \pi_A(dp|D^n), \tag{22}$$

where $\pi_A(dp|D^n) = \pi(dp|D^n, p \in A)$, and $p \in A$ is a selection rule, possibly data dependent. A rule of this kind, for example, can be averaging over several of the best models, which are indexed by $\gamma$'s having the largest marginal posteriors $\pi(\gamma|D^n)$. For example, Smith and Kohn [23] considered the use of the best model, and Sha et al. [22] averaged over the ten best models. A rule can also be defined by using the models that include the individually strongest variables. For example, include a model $\gamma$ in the posterior average if $\gamma$ includes a variable $j$ that appears more than 5 percent of time in the posterior distribution [i.e., if $\pi(\gamma_j = 1|D^n) > 0.05$]. See, for example, Lee et al. [19].

Suppose a rule $A$ has selection probability $\pi\{p \in A|D^n\} > r$ for some constant $r > 0$. Then the convergence rate of $\hat{p}_A$ can be studied by using the relations

$$d(\hat{p}_A, p^*)^2 \leq \varepsilon_n^2 + 2\pi[d(p, p^*) > \varepsilon_n | D^n]/r \tag{23}$$

and

$$P^*[d(\hat{p}_A, p^*)^2 \leq \varepsilon_n^2 + 2\delta_n/r] \geq P^*[\pi(d(p, p^*) > \varepsilon_n | D^n) \leq \delta_n], \tag{24}$$

which follows a familiar treatment based on convexity of $p \mapsto d(p, p^*)^2$ (e.g., Ghosal, Ghosh and van der Vaart [10]). The term $\delta_n/r$ can usually be taken as $e^{-r_0 n \varepsilon_n^2}$ [see result (i) of Theorem 1], which is negligible compared to $\varepsilon_n^2$ under conditions $\bar{r}_n \ln(1/\varepsilon_n^2) \prec n \varepsilon_n^2$ and $1 \leq \bar{r}_n$ of Theorem 1.

For regression purposes, a related mean estimate can be constructed as $\hat{\mu}_A(x) = \int y \hat{p}_A(y|x) \nu_y(dy)$. When binary response $y$ is considered, a classifier can be defined as $\hat{C}_A(x) = I[\hat{\mu}_A(x) > 0.5]$.

In the general case, there is no relationship bounding the $L_2$ distance $E_x(\hat{\mu}_A - \mu^*)^2$ between the estimated mean and the true mean using $d(\hat{p}_A, p^*)$, because the latter is bounded but the former is not. However, a weighted $L_2$ difference can be bounded as

$$\int \frac{(\hat{\mu}_A - \mu^*)^2}{\hat{\nu}_A + \nu^*} \nu_x(dx) \leq 2d(\hat{p}_A, p^*)^2, \tag{25}$$

where $\nu^*(x) = \int y^2 p^*(y|x)\nu_y(dy)$ and $\hat{\nu}_A(x) = \int y^2 \hat{p}_A(y|x)\nu_y(dy)$. This is obtained by noting that $(\hat{\mu}_A - \mu^*)^2 = \{\int y(\sqrt{\hat{p}_A} + \sqrt{p^*}) \cdot (\sqrt{\hat{p}_A} - \sqrt{p^*})\nu_y(dy)\}^2$ and applying the Cauchy–Schwarz inequality.

Since the denominator is at most 2 for binary $y$, the above relation actually leads to a bound for the unweighted $L_2$ distance between the means, which further leads to a bound for the classification error due to Corollary 6.2 of Devroye, Györfi and Lugosi [4]. This is summarized below and was used in proving regression and classification consistency in Jiang [15]:

$$E^*_{D^n} P^*_{(x,y)}(\hat{C}_A(x) \neq y|D^n) - L^* \tag{26}$$

$$\leq E^*_{D^n} 2\sqrt{E_x(\hat{\mu}_A - \mu^*)^2}$$

$$\leq E^*_{D^n} 2\sqrt{4d(\hat{p}_A, p^*)^2} \leq 2\sqrt{4E^*_{D^n} d(\hat{p}_A, p^*)^2} \tag{27}$$

$$\leq 4\sqrt{\varepsilon_n^2 + 2E^*_{D^n}\pi[d(p,p^*) > \varepsilon_n|D^n]/r}. \tag{28}$$

[The last step is due to (23).] Here $L^* = P^*_{(x,y)}\{C^*(x) \neq y\}$ is the "Bayes error," where $C^*(x) = I[\mu^*(x) > 1/2]$ is the ideal "Bayes rule" based on the (unknown) true mean function $\mu^*$. According to Theorem 1(iii), the term $E^*_{D^n}\pi[d(p,p^*) > \varepsilon_n|D^n]$ can be made exponentially small (of the form $e^{-c_1 n \varepsilon_n^2}$ for some $c_1 > 0$), which is negligible when compared to $\varepsilon_n^2$ as commented earlier. This implies that the error of the classification rule $\hat{C}_A(x)$ is at most $5\varepsilon_n$ above that of the optimal Bayes rule, for all large enough $n$. *So $\varepsilon_n$ also forms a rate of convergence to the optimal Bayes error for the purpose of classification.* Here the convergence rate $\varepsilon_n$ can be made to be near "finite-dimensional" (nearly $1/\sqrt{n}$) by Bayesian variable selection, despite a high dimension $K_n \sim n^\alpha \gg n$, in situations commented on earlier (e.g., Corollary 1 and Remark 1).

These convergence rate results show that even in high dimensions with $\dim(x) \gg n$, a good convergence rate can be achieved when the effect of $x$ is "sparse." For such sparse problems Bayesian variable selection can therefore help to alleviate "overfitting" or the "curse of dimensionality."

**6. Gaussian variable selection and graphical models.** In this section we will assume that $x_1^{J_n} \equiv (x_1, \ldots, x_{K_n}, y)$ are multivariate Gaussian and have been standardized to have $E(x_k) = 0$, $\text{var}(x_k) = 1$. Here $y$ is regarded as $x_{J_n}$, where $J_n = K_n + 1$. The effects of $x_{k \neq j}$ on $x_j$ are summarized by the regression coefficients $\beta_{j|k}^*$ used in the induced relation $E(x_j|x_{k \neq j}) = \sum_{k \neq j} \beta_{j|k}^* x_k$.

In Gaussian graphical models, relations among $x_1^{J_n}$ are described by a graph, such that a node corresponding to $x_j$ is only connected to a "neighborhood" $(x_k)_{k \in nb_j}$, where $nb_j$ is a subset of $\{1, \ldots, J_n\} \setminus \{j\}$, which indicates selected variables used in regression modeling of $x_j | x_{k \neq j}$. Therefore, the Bayesian variable selection technique can be used for studying the neighborhood of a variable $x_j$ (see, e.g., Dobra et al. [5]). We will consider the situation when none of the effects of $x_{k \neq j}$ on $x_j$ is exactly zero. In this case, the usual consistency of selecting the "true graph" (e.g., Meinshausen and Bühlmann [20]) will not be studied here, since the true graph is the saturated graph adopting all $K_n$ variables $x_{k \neq j}$ to explain each $x_j$. In the high-dimensional case $K_n \gg n$, such a "true model" is obviously not very useful. Nevertheless, in such a situation, Bayesian variable selection can still be shown to produce "good" models that are much simpler and yet are still "consistent," if the effects of these $K_n$ variables decay sufficiently fast (when ordered in some way). Here "consistency" is in a different sense—these simplified models, picking up only a small number out of all the $K_n$ nonzero regression coefficients, will be consistent in terms of producing probability densities "often close" to the true probability density. In this approach, one first uses Bayesian variable selection to obtain such "good" density estimates, for all $p^*(x_j|x_{k \neq j})$, $j = 1, \ldots, J_n$; then one can construct graphs to summarize the conditional independence structures corresponding to these "good" density estimates. (One can systematically decide to either include or exclude one-sided connections in these graphs (see, e.g., Meinshausen and Bühlmann [20]) when some $x_k$ is used in modeling $x_j|x_{k \neq j}$ but $x_j$ is not used in modeling $x_k|x_{j \neq k}$.)

We are interested in making inference on $J_n$ ($= K_n + 1$) (conditional) densities $p^*(x_j|x_{k \neq j})$, $j = 1, \ldots, J_n$, in order to construct a graph. We hope that the $P^*$ probability for not reliably estimating each density is small enough so that the $P^*$ probability is small for *any* density to be badly estimated. In other words, we would like to have a bound of $P^*$ probability of large errors. For now, pick any $x_j$ as the response $y$ and consider its regression on the $x_k$'s ($k \neq j$). To mimic the regression setup, we can reorder the indices of the $x_{k \neq j}$'s as $x_1^{K_n}$. We will use the prior specified in Remark 3. A result as in Theorem 1(ii), obtained when assuming uniformly bounded $|x_k|$'s, could be used for this purpose of bounding the total error out of the $J_n$ regression analyses.

In the current situation of Gaussian graphical models, however, the $x_k$'s are Gaussian instead of being uniformly bounded. In this case, for result (ii) of Theorem 1 to hold, we will change the condition on $\Delta(r_n) = \inf_{|\gamma|=r_n} \sum_{k \notin \gamma} |\beta_k^*|$ from $\Delta(r_n) \prec \varepsilon_n^2$ to $K_n \Delta(r_n) \prec \varepsilon_n^2$. This would be satisfied if $\Delta(r_n)$ decays exponentially fast in $r_n$, $r_n \succ \ln n$, and $K_n$ grows at most polynomially. After taking into account some other conditions, we obtain the following theorem.

THEOREM 3. *Consider the prior specification in Remark 3. (When selecting the neighborhood for each $x_j$, treat $x_j$ as $y$ and $x_{k \neq j}$ as $x_1^{K_n}$.) Assume that*

$$\max\left\{\sup_{\gamma: |\gamma| \leq \bar{r}_n} \mathrm{ch}_1(V_\gamma), \sup_{\gamma: |\gamma| \leq \bar{r}_n} \mathrm{ch}_1(V_\gamma^{-1})\right\} \leq B \bar{r}_n^v$$

*for some positive constants $B$ and $v$, for all large enough $\bar{r}_n$.*

*Suppose that, for each $x_j$, the effects of the other variables $x_{k \neq j}$ satisfy $\lim_{n \to \infty} \sum_{k \in \mathcal{K}_j} |\beta_{j|k}^*| < \infty$, where $\mathcal{K}_j = \{1, \ldots, K_n + 1\} \setminus \{j\}$. In addition, assume that there exists some $C' > 0$, such that, for all large enough $n$, $\inf_{\gamma \subset \mathcal{K}_j, |\gamma| = r_n} \sum_{k \in \mathcal{K}_j \setminus \gamma} |\beta_{j|k}^*| \leq e^{-C' r_n}$.*

*Assume that $n^\delta \prec K_n \prec n^\alpha$ for some $\alpha > \delta > 0$.*

*Assume also for some $\xi \in (0, 1)$*

(29) $\qquad \ln n \prec r_n \leq \bar{r}_n \prec n^b, \qquad$ *where $b < \min\{\delta, \xi, \xi/v\}$.*

*Then we have, for some constant $c'_{1,2,3} > 0$, for all sufficiently large $n$,*

(i)

$$P^*[\pi(h_j \leq n^{-(1-\xi)/2} | D^n) \geq 1 - e^{-c'_1 n^\xi}, j = 1, \ldots, K_n + 1] \geq 1 - n^\alpha e^{-c'_2 n^\xi}$$

*and*

(ii)

$$P^*[\hat{h}_{j, A_j} \leq c'_3 n^{-(1-\xi)/2}, j = 1, \ldots, K_n + 1] \geq 1 - n^\alpha e^{-c'_2 n^\xi}.$$

*Here we define, for $j = 1, \ldots, K_n + 1$,*

(30) $\quad h_j = \left\{\int_{\Re^{K_n+1}} |p(x_j | x_{k \neq j})^{1/2} - p^*(x_j | x_{k \neq j})^{1/2}|^2 p^*(x_{k \neq j}) \, dx_1^{K_n+1}\right\}^{1/2},$

(31)
$$\hat{h}_{j, A_j} = \left\{\int_{\Re^{K_n+1}} |\hat{p}_{A_j}(x_j | x_{k \neq j})^{1/2} - p^*(x_j | x_{k \neq j})^{1/2}|^2 \right.$$
$$\left. \times p^*(x_{k \neq j}) \, dx_1^{K_n+1}\right\}^{1/2},$$

*where $p^*$ represents the true density and $\hat{p}_{A_j}$ is a selected posterior estimate [as defined in (22)] corresponding to a selection rule $A_j$, such that the selection probability $\pi(p \in A_j | D^n) > r$ for some $r > 0$.*

Therefore, a near finite-dimensional rate of convergence can be achieved (for some small $\xi > 0$), jointly for all neighborhoods of $x_j$, $j = 1, \ldots, K_n + 1$, despite the fact that $K_n$ can follow a large power of $n$.

**7. General prior.** In this section we consider the case $|x_j| \leq 1$ for all $j$ and mainly focus on the GLM models as described in Section 2, where $a(h)$ and $b(h)$ contain no additional parameters other than $h$. (Similar conditions and results can be formulated for normal linear regression with unknown error variance.)

Here we consider the general conditions on the prior $\pi(\gamma, \beta_\gamma)$ for producing rate of convergence $\varepsilon_n$, which is a sequence in $n$, which we assume to satisfy $\varepsilon_n \in (0, 1]$ for Conditions (N) and (O) below.

Condition (N) requires a not too little prior to be placed over a very small neighborhood of the true density $p^*$. Condition (O) requires a very little prior to be placed outside of a region that is not too complex in some sense.

CONDITION (N) (*For prior $\pi$ on an approximation neighborhood*). Assume that a sequence of (nonempty) models $\gamma_n$ exists such that, as $n$ increases,

$$\sum_{j \notin \gamma_n} |\beta_j^*| \prec \varepsilon_n^2, \tag{32}$$

and for any sufficiently small $\eta > 0$, there exists $N_\eta$ such that, for all $n > N_\eta$, we have

$$\pi(\gamma = \gamma_n) \geq e^{-n\varepsilon_n^2/8} \tag{33}$$

and

$$\pi(\beta_\gamma \in M(\gamma_n, \eta) | \gamma = \gamma_n) \geq e^{-n\varepsilon_n^2/8}, \tag{34}$$

where $M(\gamma_n, \eta) = (\beta_j^* \pm \eta \varepsilon_n^2 / |\gamma_n|)_{j \in \gamma_n}$.

CONDITION (O) (*For prior $\pi$ outside of a not-too-complex region*). Let $D(R) = 1 + R \sup_{|h| \leq R} |a'(h)| \cdot \sup_{|h| \leq R} |\psi(h)|$ for any $R > 0$. There exist some $C_n > 0$ and some $\bar{r}_n$ satisfying $1 \leq \bar{r}_n < K_n$, such that

$$\bar{r}_n \ln(1/\varepsilon_n^2) \prec n\varepsilon_n^2, \tag{35}$$

$$\bar{r}_n \ln K_n \prec n\varepsilon_n^2, \tag{36}$$

$$\bar{r}_n \ln D(\bar{r}_n C_n) \prec n\varepsilon_n^2. \tag{37}$$

Furthermore, for all large enough $n$, the following two equations hold:

$$\pi(|\gamma| > \bar{r}_n) \leq e^{-4n\varepsilon_n^2}, \tag{38}$$

and for all $\gamma$ such that $|\gamma| \leq \bar{r}_n$, for all $j \in \gamma$,

$$\pi(|\beta_j| > C_n | \gamma) \leq e^{-4n\varepsilon_n^2}. \tag{39}$$

These conditions allow a larger variety of priors. For example, one can use a uniform prior of $\gamma$ over all models with complexity $|\gamma| \leq \bar{r}_n$, where $\bar{r}_n$ can be taken to follow some rate of growth depending on the convergence rate $\varepsilon_n$ desired, and depending on the "bias" rate $\Delta(r) = \inf_{\gamma: |\gamma|=r} \sum_{j \notin \gamma} |\beta_j^*|$.

Before we truncated $\pi(\gamma)$ such that $\pi[|\gamma| > \bar{r}_n] = 0$. This may not be desirable since we are forbidding the model to be too complex in the prior. We here notice that this truncation is not necessary. We can allow the prior to propose very complicated models with large $|\gamma|$, as long as the prior probability of $|\gamma| > \bar{r}_n$ is sufficiently small.

THEOREM 4 (Convergence rate under general prior). *For GLM models with bounded covariates $|x_j| \leq 1$ for all $j$, suppose the true regression coefficients satisfy $\lim_{n \to \infty} \sum_{j=1}^{K_n} |\beta_j^*| < \infty$.*

*Let $\varepsilon_n \in (0, 1]$ be a sequence such that $n\varepsilon_n^2 \to \infty$. Denote $d(p, p^*)^2 = \int \int |p(y, x|\gamma, \beta_\gamma)^{1/2} - p^*(y, x)^{1/2}|^2 \nu_y(dy) \nu_x(dx)$. If the prior specification satisfies both Conditions* (N) *and* (O), *then we have the following (successively stronger) results:*

(i)
$$\lim_{n \to \infty} P^*\{\pi_n[d(p, p^*) \leq 4\varepsilon_n | D^n] \geq 1 - 2e^{-n\varepsilon_n^2/4}\} = 1.$$

(ii) *For all sufficiently large $n$,*
$$P^*\{\pi[d(p, p^*) > 4\varepsilon_n | D^n] \geq 2e^{-n\varepsilon_n^2/4}\} \leq 2e^{-n\varepsilon_n^2/4}.$$

(iii) *For all sufficiently large $n$,*
$$E_{D^n}^* \pi[d(p, p^*) > 4\varepsilon_n | D^n] \leq 4e^{-n\varepsilon_n^2/2}.$$

Results (i), (ii) and (iii) of this theorem will be proved by verifying some sufficient conditions for posterior convergence (to be summarized at the beginning of Section 8). These results, with bounded $x_j$'s, will then be used to prove all the previous results on convergence rates, when specific priors as given in Section 3 are used. The only exception is the result in Section 6, where $x_j$'s are jointly normal; they will be obtained by directly verifying the sufficient conditions in Section 8.

These conditions below are based on the Hellinger metric entropy and will be used to obtain posterior convergence rates under the GLM framework. Note that the method involved here is different from that in Jiang [15], who uses the Hellinger bracketing entropy and its upper bound of a parametric covering number (see, e.g., Theorem 3, Lee [18]). That method does not directly apply to modeling unbounded responses such as Gaussian and Poisson responses. (When applied to, e.g., Poisson regression, the upper bound of the bracketing entropy would require a too small restricted parameter space, on which the prior would place a nonnegligible probability.)

**8. Proofs.** We first use a proposition to summarize a set of sufficient conditions for establishing rates of posterior convergence. These serve as just one possible set of working conditions that we find convenient to use here, through which we have established our results; there exist several other alternatives, possibly with more relaxed conditions, for example, in Ghosal, Ghosh and van der Vaart [10] or Zhang [26].

Suppose $\mathcal{P}_n$ is a sequence of sets of probability densities. (For each $n$, denote $\mathcal{P}_n^c$ as the complement—the set of densities not in $\mathcal{P}_n$.) Suppose $\varepsilon_n$ is a sequence of positive numbers.

Suppose $N(\varepsilon_n, \mathcal{P}_n)$ is the minimal number of Hellinger balls of radius $\varepsilon_n$ that are needed to cover $\mathcal{P}_n$. [I.e., $N(\varepsilon_n, \mathcal{P}_n)$ is the minimum of all $k$ such that there exist $S_j = \{p : d(p, p_j) \leq \varepsilon_n\}$, $j = 1, \ldots, k$, such that $\bigcup_{j=1}^k S_j \supset \mathcal{P}_n$, where $d(p,q) = \sqrt{\int (\sqrt{p} - \sqrt{q})^2}$ denotes the Hellinger distance between densities $p$ and $q$.]

Let the components of $D^n = (w^{(1)}, \ldots, w^{(n)})$ be i.i.d. with true density $p^*$, where $\dim(w^{(1)})$ and $p^*$ can depend on $n$. Denote $\pi(\cdot)$ as the prior (which is allowed to depend on $n$ by using, e.g., an increasing number of parameters to parameterize the density as $n$ increases), $\pi(\cdot|D^n)$ as the posterior and $\hat{\pi}(\varepsilon) = \pi[d(p, p^*) > \varepsilon | D^n]$ for each $\varepsilon > 0$. Define the KL difference as $d_0(p, p^*) = \int p^* \ln(p^*/p)$. Define also a $d_t$ difference as $d_t(p, p^*) = t^{-1}(\int p^*(p^*/p)^t - 1)$ for any $t > 0$, which is used in, for example, Wong and Shen [25]. (Note that $d_t$ decreases to $d_0$ as $t$ decreases toward 0.)

Denote $P^*$ and $E^*$ as the respective probability measure and the expectation for the data $D^n$.

Define the following conditions:

(a) $\ln N(\varepsilon_n, \mathcal{P}_n) \leq n\varepsilon_n^2$ for all sufficiently large $n$;
(b) $\pi(\mathcal{P}_n^c) \leq e^{-2n\varepsilon_n^2}$ for all sufficiently large $n$;
(c) for all small enough $\gamma > 0$ and $r > 0$, there exists $N_{\gamma,r}$ such that for all $n \geq N_{\gamma,r}$, $\pi[p : d_0(p, p^*) \leq \gamma \varepsilon_n^2] \geq e^{-rn\varepsilon_n^2}$;
(d) $\pi[p : d_t(p, p^*) \leq \varepsilon_n^2/4] \geq e^{-n\varepsilon_n^2/4}$ for all sufficiently large $n$, for some $t > 0$.

PROPOSITION 1. *Suppose $n\varepsilon_n^2 \succ 1$. Then, under (a), (b) and (c), we have:*

(i)
$$\lim_{n \to \infty} P^*[\hat{\pi}(4\varepsilon_n) < 2e^{-n\varepsilon_n^2 \min\{1, b/2\}}] = 1;$$

*under (a), (b) and (d) (for some $t > 0$), we have*
(ii)
$$P^*[\hat{\pi}(4\varepsilon_n) \geq 2e^{-n\varepsilon_n^2 \min\{1/2, t/4\}}] \leq 2e^{-n\varepsilon_n^2 \min\{1/2, t/4\}},$$

(iii)

$$E^*\hat{\pi}(4\varepsilon_n) \leq 4e^{-n\varepsilon_n^2 \min\{1, t/2\}}.$$

The proof of this proposition follows the spirit of Ghosal, Ghosh and van der Vaart [10]. The details are omitted here and are included in a technical report (Jiang [14]).

PROOF OF THEOREM 4. We prove result (iii) only, since it implies (ii) by Markov's inequality, which further implies (i). Result (iii) is proven by applying Proposition 1 with $t = 1$. The proof is completed by checking conditions (d), (a) and (b) below.

*Checking condition* (d) *for $t = 1$.* Denote the GLM density as $f(y, h) = \exp\{a(h)y + b(h) + c(y)\}$. Then $p^* = f(y, h^*)$, where $h^* = x^T \beta^* = \sum_{j=1}^{K_n} x_j \beta_j^*$. Let $p_\gamma = f(y, h_\gamma)$, where $h_\gamma = x_\gamma^T \beta_\gamma = \sum_{j \in \gamma_n} x_j \beta_j$, where $\gamma_n$ is the model in Condition (N).

When $h^*$ and $h_\gamma$ are close enough, $d_t(p_\gamma, p^*)$ (for $t = 1$) can be put in a form $d_t(p_\gamma, p^*) = E_x g(h^i)(h^* - h_\gamma)$, by integrating out $y$ and applying a first-order Taylor expansion. Here $g$ is a continuous derivative function in a neighborhood of $h^*$ and $h^i$ is an intermediate point between $h^*$ and $h_\gamma$. Note that $|h^i - h^*| \leq |h_\gamma - h^*| \leq |\sum_{j \notin \gamma_n} x_j \beta_j^*| + |\sum_{j \in \gamma_n} x_j (\beta_j - \beta_j^*)| \leq \Delta_n + r_n \delta_n$, when the $x_j$'s are bounded by 1 and $\beta_j \in (\beta_j^* \pm \delta_n)$ for all $j \in \gamma_n$. Here $r_n = |\gamma_n|$ and $\Delta_n = \sum_{j \notin \gamma_n} |\beta_j^*|$ is assumed to satisfy $\Delta_n \prec \varepsilon_n^2$.

For sufficiently small $r_n \delta_n$, $|g(h^i)|$ is bounded since $|h^i| \leq |h^*| + |h^i - h^*| \leq B_0 + \Delta_n + r_n \delta_n$ is bounded, where $B_0 = \lim_{n \to \infty} \sum_{j=1}^{K_n} |\beta_j^*|$. Then $d_t(p_\gamma, p^*) \leq C(\Delta_n + r_n \delta_n)$ for some constant $C$, for all small enough $r_n \delta_n$.

We will take $\delta_n = \eta \varepsilon_n^2 / |\gamma_n|$ for some small enough $\eta > 0$. This will make $d_t \leq \varepsilon_n^2 / 4$ for all large enough $n$, since $\Delta_n \prec \varepsilon_n^2$.

This implies that the set of densities $S = \{p(\cdot|\gamma_n, \beta) : \beta \in (\beta_j^* \pm \delta_n)_{j \in \gamma_n}\}$ is contained in $T = \{p : d_t(p, p^*) \leq \varepsilon_n^2 / 4\}$. The conditions on $\pi(\gamma_n)$ and $\pi(\beta \in (\beta_j^* \pm \delta_n)_{j \in \gamma_n} | \gamma_n)$ then imply that $\pi(T) \geq \pi(S) \geq e^{-n\varepsilon_n^2/4}$ for all large $n$, confirming condition (d).

*Checking condition* (a). Each density $p$ is labeled by a model index $\gamma$ and the corresponding regression coefficients $\beta_j$. We will define $\mathcal{P}_n$ as the set of densities that can be represented with $|\gamma|$ (the number of nonzero regression parameters) being at most $\bar{r}_n$, and with each parameter $|\beta_j| \leq C_n$.

The corresponding space of regression parameters can be covered by small $\ell_\infty$ balls of the form $B = (v_j \pm \delta)_{j=1}^{K_n}$, of radius $\delta > 0$. For each model $\gamma$ in $\mathcal{P}_n$, there are $|\gamma|$ nonzero components of $\beta_j$, valued in $\pm C_n$. It takes at most $[(2C_n)/(2\delta) + 1]^{|\gamma|}$ balls to cover the parameter space of model $\gamma$ in $\mathcal{P}_n$. [The centers of these balls can be taken inside the parameter space of model $\gamma$, so that each center $v = (v_j)_1^{K_n}$ has components satisfying $v_j = 0$ $\forall j \notin \gamma$ and $|v_j| \leq C_n$ $\forall j \in \gamma$.]

There are at most $K_n^r$ models of size $|\gamma| = r$, and $r = 0, 1, 2, \ldots, \bar{r}_n$. These show that $N(\delta)$, the number of size-$\delta$ balls needed to cover the space of regression parameters for $\mathcal{P}_n$, is at most $\sum_{r=0}^{\bar{r}_n} K_n^r [(2C_n)/(2\delta) + 1]^r$, which is bounded above by $(\bar{r}_n + 1)(K_n(C_n/\delta + 1))^{\bar{r}_n}$.

Given any density in $\mathcal{P}_n$, it can be represented by a set of regression parameters $(u_j)_1^{K_n}$ falling in one of these $N(\delta)$ balls, say, ball $B = (v_j \pm \delta)_{j=1}^{K_n}$, where $u_j$ and $v_j$ are zero for the same set of components $\gamma$, where $|\gamma| \leq \bar{r}_n$.

Consider the corresponding GLM densities $p_{u,v} = \exp\{ya(h_{u,v}) + b(h_{u,v}) + c(y)\}$, where $h_u = \sum_{j=1}^{K_n} u_j x_j$ and $h_v = \sum_{j=1}^{K_n} v_j x_j$. Then the Hellinger distance $d(p_u, p_v) \leq \{d_0(p_u, p_v)\}^{1/2}$, where the KL difference $d_0(p_u, p_v) = E_x \int p_v (\ln p_v - \ln p_u) \nu_y(dy)$. After integration in $y$, one can apply a Taylor expansion and show that $d_0(p_u, p_v) \leq E_x(a'(h^i)\psi(h_v) + b'(h^i))(h_v - h_u)$, where $\psi = -b'/a'$ and $h^i$ is an intermediate point between $h_v$ and $h_u$. Note that $u$ and $v$ both have components bounded in value by $C_n$ and they have zero components out of a same set, say, $\gamma$, such that $|\gamma| \leq \bar{r}_n$. Therefore, $h_v$ and $h_u$ (and therefore, also $h^i$) are bounded above by $\bar{r}_n C_n$. Note also that $|h_v - h_u| = |\sum_{j \in \gamma} x_j(v_j - u_j)| \leq \bar{r}_n \delta$, since $|x_j| \leq 1$, $|v_j - u_j| \leq \delta$ and $|\gamma| \leq \bar{r}_n$. Therefore,

$$(40) \qquad d_0(p_u, p_v) \leq 2 \sup_{|h| \leq \bar{r}_n C_n} |a'(h)| \sup_{|h| \leq \bar{r}_n C_n} |\psi(h)| \bar{r}_n \delta$$

and

$$(41) \qquad d(p_u, p_v) \leq \left\{ 2 \sup_{|h| \leq \bar{r}_n C_n} |a'(h)| \sup_{|h| \leq \bar{r}_n C_n} |\psi(h)| \bar{r}_n \delta \right\}^{1/2}.$$

So $d(p_u, p_v) \leq \varepsilon_n$ if $\delta = \varepsilon_n^2 / \{2 \sup_{|h| \leq \bar{r}_n C_n} |a'(h)| \sup_{|h| \leq \bar{r}_n C_n} |\psi(h)| \bar{r}_n \}$. Therefore, density $p_u$ in $\mathcal{P}_n$ falls in a Hellinger ball of size $\varepsilon_n$, centered at $p_v$. There are at most $N(\delta)$ such balls, because each center $p_v$ is the density corresponding to the parameter $v$, which is the center of $B$, one of the at most $N(\delta)$ balls used to cover the restricted parameter space.

Therefore, the Hellinger covering number

$N(\varepsilon_n, \mathcal{P}_n) \leq N(\delta)$

$$(42) \qquad \leq (\bar{r}_n + 1) K_n^{\bar{r}_n} \left( 1 + 2\varepsilon_n^{-2} \sup_{|h| \leq \bar{r}_n C_n} |a'(h)| \sup_{|h| \leq \bar{r}_n C_n} |\psi(h)| \bar{r}_n C_n \right)^{\bar{r}_n}$$

$$\leq (2K_n^2 D(\bar{r}_n C_n)/\varepsilon_n^2)^{\bar{r}_n},$$

if $0 < \varepsilon_n \leq 1$ and $1 \leq \bar{r}_n < K_n$. Therefore, the conditions in Condition (O) guarantee that $\ln N(\varepsilon_n, \mathcal{P}_n) \prec n\varepsilon_n^2$ for all large enough $n$, proving condition (a).

*Checking condition* (b). For the $\mathcal{P}_n$ defined above, the prior on the complement $\pi(\mathcal{P}_n^c) \leq \pi[|\gamma| > \bar{r}_n] + \sum_{\gamma:|\gamma|\leq \bar{r}_n} \pi[\gamma]\pi(\bigcup_{j\in\gamma}[|\beta_j| > C_n]|\gamma)$, which is at most $\pi[|\gamma| > \bar{r}_n] + \max_{\gamma:|\gamma|\leq \bar{r}_n} \pi(\bigcup_{j\in\gamma}[|\beta_j| > C_n]|\gamma)$. This is, due to Condition (O), at most $(1+\bar{r}_n)e^{-4n\varepsilon_n^2} = e^{\ln(1+\bar{r}_n)-4n\varepsilon_n^2} \leq \exp(-4n\varepsilon_n^2/2)$ for all large enough $n$. Here we have used $1 \leq \bar{r}_n < K_n$, so that $\ln(1+\bar{r}_n) \leq \bar{r}_n \ln K_n \prec n\varepsilon_n^2$ due to Condition (O). This proves condition (b). □

PROOF OF THEOREM 1. We apply Theorem 4 with $\varepsilon_n$ replaced by $\varepsilon_n'$, so that the Hellinger neighborhood will take a radius $4\varepsilon_n'$. This can be later rescaled to obtain the results in Theorem 1 concerning a radius $\varepsilon_n$, by setting $\varepsilon_n = 4\varepsilon_n'$ or $\varepsilon_n' = \varepsilon_n/4$.

For Condition (O): with the prior in Section 3, the condition on $\pi[|\gamma| > \bar{r}_n]$ is trivially satisfied, since it is zero due to truncation. We will take $C_n = \sqrt{\tilde{B}_n n\varepsilon_n^2}$ so that the condition on $\bar{r}_n \ln D(\bar{r}_n C_n)$ is satisfied. The condition on $\pi[|\gamma_j| > C_n|\gamma]$ is checked by using Mills' ratio. It is at most $2e^{-C_n^2/(2\tilde{B}_n)}/\sqrt{2\pi C_n^2/\tilde{B}_n}$, which is therefore less than $e^{-n\varepsilon_n^2/4} = e^{-4n(\varepsilon_n')^2}$ for all large enough $n$, as required by Condition (O). Here $\tilde{B}_n$ is an upper bound on the prior variance of $\beta_j$ under model $\gamma$ with $|\gamma| \leq \bar{r}_n$, and $n\varepsilon_n^2 \succ 1$. All other conditions in Condition (O) are satisfied.

For Condition (N): Take the sequence of models $\gamma_n$ such that, for each $n$, $\gamma = \gamma_n$ reaches its infimum in $\Delta(r_n) = \inf_{\gamma:|\gamma|=r_n} \sum_{j:j\notin\gamma} |\beta_j^*|$. Then $\sum_{j\notin\gamma_n} |\beta_j^*| = \Delta(r_n) \prec n\varepsilon_n^2$.

For the condition on the prior $\pi[\beta \in (\beta_j^* \pm \eta\varepsilon_n^2/r_n)_{j\in\gamma_n}|\gamma_n]$, use the normality of the prior and obtain the lower bound $|2\pi V_{\gamma_n}|^{-1/2}e^{-0.5\beta^T V_{\gamma_n}^{-1}\beta}(\eta\varepsilon_n^2/r_n)^{r_n}$ for some intermediate value $\beta$ achieving the infimum of the density over $(\beta_j^* \pm \eta\varepsilon_n^2/r_n)_{j\in\gamma_n}$.

Note that $\beta^T V_{\gamma_n}^{-1}\beta \leq \|\beta\|^2 B(r_n) \leq (\sum_{j\in\gamma_n}|\beta_j|)^2 B(r_n) \leq C_1 B(r_n)$ for some constant $C_1 > 0$, since the eigenvalues of $V_{\gamma_n}^{-1}$ are at most $B(r_n)$ (for all large enough $n$), and the Euclidean norm $\|\beta\| \leq \sum_{j\in\gamma_n}|\beta_j| \leq \lim_{n\to\infty}\sum_{j=1}^{K_n}|\beta_j^*| + r_n\eta\varepsilon_n^2/r_n$ is bounded. Note also that $|2\pi V_{\gamma_n}|^{-1/2} \geq e^{-C_2 r_n - C_3 r_n \ln \bar{B}(r_n)}$ for some constant $C_2$ and some constant $C_3 > 0$, due to the eigenvalues of $V_{\gamma_n}$ being bounded above by $\bar{B}(r_n)$ (for all large enough $n$).

Therefore,

$$
(43) \quad \begin{aligned} &\pi[\beta \in (\beta_j^* \pm \eta\varepsilon_n^2/r_n)_{j\in\gamma_n}|\gamma_n] \\ &\geq \exp\{-C_2 r_n - C_3 r_n \ln \bar{B}(r_n) - 0.5 C_1 B(r_n) - r_n \ln(r_n/(\eta\varepsilon_n^2))\}. \end{aligned}
$$

This will be greater in order than any $e^{-cn\varepsilon_n^2}$ ($c > 0$), satisfying a requirement of Condition (N), since $r_n$, $r_n \ln \bar{B}(r_n)$ and $B(r_n)$ are all smaller than $n\varepsilon_n^2$ in order, and so are $r_n \ln r_n \leq \bar{r}_n \ln K_n$ and $r_n \ln(1/\varepsilon_n^2) \leq \bar{r}_n \ln(1/\varepsilon_n^2)$.

Now consider the condition on $\pi(\gamma_n)$. Note that the $\gamma_n$ chosen is such that $|\gamma_n| = r_n$, where $r_n \ (\leq \bar{r}_n)$ is the expected size of the model $\tilde{\gamma} = \tilde{\gamma}_1^{K_n}$ proposed by the prior *before* truncation. The prior specification of $\gamma_1^{K_n}$ (in Section 3) is i.i.d. binary with $\pi(\tilde{\gamma}_j = 1) = r_n/K_n$. For the condition on $\pi(\gamma = \gamma_n)$ to hold, it suffices for us to show that (*) for any $c > 0$, $\pi(\tilde{\gamma} = \gamma_n) > e^{-cn\varepsilon_n^2}$ for all large enough $n$. This is because $\pi(\gamma = \gamma_n)$ cannot be smaller, since it is obtained by truncation of $\tilde{\gamma}$, and truncation increases the probability of all allowed configurations (note that $|\gamma_n| \leq \bar{r}_n$).

Now $|\gamma_n| = r_n \leq \bar{r}_n$ implies that there are $r_n$ out of $K_n$ $\tilde{\gamma}_j$'s equal to 1, with the rest being 0. The probability is therefore $\pi(\tilde{\gamma} = \gamma_n) = (r_n/K_n)^{r_n}(1 - r_n/K_n)^{K_n - r_n}$. Since $r_n/K_n \prec 1$, we have $\ln \pi(\tilde{\gamma} = \gamma_n) \sim r_n \ln(r_n/K_n) \geq -r_n \ln K_n \ (r_n \geq 1)$, where $r_n \ln K_n \prec n\varepsilon_n^2$. This leads to claim (*). □

PROOF OF THEOREM 3. It suffices for us to prove (**) result (ii) of Proposition 1 in a regression setup for normal dispersion models with $K_n$ Gaussian covariates. We will take $\varepsilon_n \sim n^{-(1-\xi)/2}$ for some $\xi \in (0, 1)$. Then result (i) of Theorem 3 can be obtained by a union bound over $K_n + 1$ regressions, treating each of the $K_n + 1$ $x_j$'s in turn as the response $y$, and $x_{k \neq j}$ as $x_1^{K_n}$. Result (ii) of Theorem 3 can be obtained by using bounds of the form (24).

We prove (**) by directly applying Proposition 1 and verifying conditions (a), (b) and (d) (for $t = 1$). The details are omitted here and are included in a technical report (Jiang [14]) in order to save space. □

**9. Discussion.** Bayesian variable selection (BVS) handles high-dimensional regression by using a suitable prior to propose lower-dimensional models which select a few explanatory variables out of the many ($K_n$) candidates. For generalized linear models, we have shown that (see, e.g., Remark 1) a near finite-dimensional convergence rate $\varepsilon_n$ can be obtained, even when the number of candidate variables $K_n$ grows as any high power $n^\alpha$ of the sample size $n$. Such a good rate $\varepsilon_n$ is derived assuming an exponentially decaying tail $\Delta(r_n) = \inf_{\gamma: |\gamma| = r_n} \sum_{j \notin \gamma} |\beta_j^*|$. This includes as a special case the situation when only a fixed and finite number of true regression coefficients ($\beta_j^*$'s) are nonzero. On the other hand, it also allows more realistic situations with many small $|\beta_j^*|$'s, none of which is exactly zero. The rates we obtain here are infinitesimally weaker than the finite-dimensional rate $n^{-1/2}$. We suspect that the exact rate $n^{-1/2}$ cannot be achieved in the setup that we consider, since the priors we use need to propose models of dimension $r_n$ increasing to infinity as $n$ increases (even though $r_n \ll n$). This is for the purpose of being able to approximate a true model to any precision. With such increasing model dimensions, we suspect that the exact $n^{-1/2}$-rate cannot be achieved in any way.

Although we have only considered in detail the situation with an exponentially decaying $\Delta(\cdot)$, the more general framework of, for example, Theorem 2 allows us to treat other situations of $\Delta(\cdot)$ as well. For example, when $\Delta(\cdot)$ follows an inverse power law, the convergence rate $\varepsilon_n$ can be somewhat slower. However, even in such situations, BVS can still exhibit some "resistance against overfitting" when $K_n$ is large. Not only can posterior consistency be still achieved when $\lim_{n\to\infty} \sum_{j=1}^{K_n} |\beta_j^*| < \infty$, but also the convergence rates will not be directly linked to the large dimension $K_n$—they will be related to the sizes of the $|\beta_j^*|$'s instead.

An Associate Editor raised the interesting question whether the sparseness conditions for the true regression coefficients can be extended to a form of $\ell_k$-summability for $k > 1$ (such as $\ell_2$). We do not have a general answer, except in an analytically-friendly special case as follows: The true model is $y \sim N(x^T \beta^*, 1)$ (it can be extended to allow a dispersion parameter), such that $Exx^T$ forms an identity matrix, or more generally, $Exx^T$ and its inverse both have bounded eigenvalues. The prior proposes fitted models of the form $y \sim N(x_\gamma^T \beta_\gamma, 1)$, according to the prior specification in Section 3. For this example, by a treatment parallel to the current paper, we can accommodate $\beta^*$ that is $\ell_2$-summable but not $\ell_1$-summable, such as $\beta^* = (j^{-1})_1^{K_n}$, resulting in a possibly slower rate for posterior convergence. On the other hand, when $\beta^* = (j^{-1/2})_1^{K_n}$, which is $\ell_k$-summable for $k > 2$ but not $\ell_2$-summable, the current approach does not work. Roughly speaking, we would need to use very complicated fitted models of size $|\gamma| \sim K_n$ to approximate the true density, in order to obtain a nonzero prior probability over a small neighborhood of the true model. Then the complexity/entropy conditions [e.g., equation (10)], which imply $|\gamma| \ln K_n \prec n$, could not be satisfied for such fitted models of size $|\gamma| \sim K_n$ in a high-dimensional setting $K_n > n$.

Although the topic of our paper is Bayesian, it is noted that the use of $\ell_1$-type conditions here is related to some other work in the frequentist approach. Our paper is closer to Bühlmann [2] in the sense that both assume a true model satisfying some $\ell_1$ summability condition, while the fitted model (boosting for Bühlmann [2] and BVS for the current paper) does not use an $\ell_1$ constraint or penalization. The fitted models in this paper are proposed according to a prior that uses i.i.d. binary distributions (with a small selection probability) when selecting the candidate variables. This may be regarded as a nondeterministic way of penalizing the $\ell_0$ norm of $\beta$ (or the number of nonzero regression coefficients) of the fitted models. On the other hand, in Greenshtein and Ritov [12], Greenshtein [11] and Meinshausen and Bühlmann [20], the fitted models (instead of the true models) are subject to an $\ell_1$ constraint or penalization. In the more general framework of persistence in Greenshtein and Ritov [12] and Greenshtein [11], the true models actually do not need to satisfy the $\ell_1$ summability condition.

The current paper focuses on fitting a density with Bayesian variable selection (BVS). A referee raised some interesting questions about the use of BVS when the main goal is selecting the variables. In some sense the current paper does prove that the method of BVS will provide "good" sets of variables, based on which good predictive performance, for example, in classification, can be achieved; see discussion in Section 5. The paper focused on the more realistic situation when there is no simple true model with many zero regression coefficients. All variables may have *some* effects, more or less. So the problem is not to select a "true" model (which would be the full model) but a "good" model (possibly much simpler than the full model) that achieves good performance for prediction, regression or density estimation. In this sense the paper does address variable selection and shows that BVS provides "good" sets of variables with high posterior probability. What will happen when there *does* exist a small true model, for example, when some regression coefficients are bounded away from zero, while the rest are exactly zero? We conjecture that, with high probability, BVS will select all the "relevant" variables with nonzero regression coefficients, but it may also include some "irrelevant" variables, with small regression coefficients proposed by the posterior. A truncation scheme similar to thresholding may be used to screen out the "irrelevant" variables, if necessary. However, we leave this as an open question, since such a scenario, being more idealized but still very interesting, is not within the main scope of the current paper.

Another future work may be to consider the (generalized) linear structure of the fitted models in a misspecified framework such as in Kleijn and van der Vaart [16], so that the true model may be nonlinear. On the other hand, one should note that nonlinearity may be treated even under the linear framework of the true model. This can be done by including higher order terms, interactions, regression spline terms with various knot-locations, and so on (see, e.g., Smith and Kohn [23]).

**Acknowledgments.** I wish to thank Professor Eitan Greenshtein for providing a related technical report. I also wish to thank the reviewers and the Editors for their insightful questions and comments.

DEPARTMENT OF STATISTICS
NORTHWESTERN UNIVERSITY
EVANSTON, ILLINOIS 60208
USA
E-MAIL: wjiang@northwestern.edu